\documentclass[a4paper,11pt]{article}

\usepackage[utf8x]{inputenc}
\usepackage{amsmath,amsfonts,amsthm,amssymb}
\usepackage{fullpage}
\usepackage{color}

\providecommand{\R}{\mathbb{R}}
\providecommand{\dx}{\, \mathrm{d} x}

\def\ds{\displaystyle}
\def\eps{\varepsilon}
\newtheorem{theorem}{Theorem}
\newtheorem{lemma}{Lemma}
\def\les{\lesssim}
\def\ges{\gtrsim}
\title{Nucleation barriers at corners for cubic-to-tetragonal phase transformation}

\author{Peter Bella
\and
        Michael Goldman\footnote{Max Planck Institute for Mathematics in the Sciences, Leipzig (Germany), email: bella@mis.mpg.de,  goldman@mis.mpg.de}
        }

\begin{document}

\maketitle

\begin{abstract}
 We are interested in the energetic cost of a martensitic inclusion of volume $V$ in austenite for the cubic-to-tetragonal phase transformation. In contrast with the work of [Kn\"upfer, Kohn, Otto: Comm. Pure Appl. Math. 66 (2013), no. 6, 867--904], we consider domain with a corner and obtain a better scaling law for the minimal energy ($E_{min} \sim \min(V^{2/3},V^{7/9})$). Our predictions are in a good agreement with physical experiments where nucleation of martensite is usually observed near the corners of the specimen. 
\end{abstract}

\section{Introduction}
In this short note we study the energetic cost of a nucleation for the cubic-to-tetragonal phase transformation in domains with a corner. We prove that in the geometrically linear setting, it is energetically more favorable to start the nucleation at the corner of the domain compared to the nucleation inside the specimen. This is in good agreement with experiments and with related results for cubic-to-orthorhombic phase transitions \cite{balKoumaSei}.

We  work in the framework of geometrically linear elasticity so that the elastic energy of our material depends only on the symmetric part $e(u):=\frac{1}{2}(\nabla^t u+\nabla u)$  of the gradient of the displacement $u : \Omega \to \R^3$ from a reference configuration (see~\cite{KKO}). We choose the reference lattice to be the austenite phase so that the stress-free strains are given by $e^{(0)}=0$ in the austenite and by $e^{(1)}$, $e^{(2)}$, and $e^{(3)}$ in the three variants of the martensite (see Section \ref{secmod} for the definition of $e^{(i)}$). We  introduce the characteristic functions of the martensitic variants
\begin{equation}\label{chiBV}
\chi_i \in BV(\Omega, \{ 0, 1 \}), \textrm{ and } \chi_1 + \chi_2 + \chi_3 \le 1.
\end{equation}

The volume of the martensitic inclusion is then given by
\begin{equation}\label{volume}
 V := \int_\Omega \chi_1 + \chi_2 + \chi_3 \dx.
\end{equation}
The normalized elastic energy is given by: 
\begin{equation*}
 E_{elast}[\chi] = \inf_{u \in H^1(\Omega,\R^3)} \int_\Omega \left| e(u) - \sum_{i=1}^3 \chi_i e^{(i)} \right|^2 \dx.
\end{equation*}
We also include an interfacial energy which penalizes the creation of interfaces between austenite and martensites or between two variants of martensites (see \cite{bhat}):
\begin{equation*}
 E_{interf}[\chi] := \int_\Omega \left( |\nabla \chi_1| + |\nabla \chi_2| + |\nabla \chi_3| \right) \dx,
\end{equation*}
where for $\chi \in BV(\Omega, \{ 0, 1 \})$, $\int_{\Omega} |\nabla \chi|$ denotes the total variation of $\chi$.
The total energy is then given by
\begin{equation*}
 E[\chi] := E_{interf}[\chi] + E_{elast}[\chi].
\end{equation*}
We are interested in the energetic cost of the formation of a nucleus of martensite of volume $V$ inside the austenite. This energetic cost determines the energy barrier which needs to be overcome to start the nucleation of martensite in the austenite. In \cite{KKO}, it was proven that when $\Omega=\R^3$, 
\[\inf_{\chi \textrm{ satisfies } \eqref{chiBV},\eqref{volume}} E[\chi] \sim \left\{ \begin{array}{ll} V^{2/3} & \textrm{ if } V \le 1 \\ V^{9/11} & \textrm{ if } V \ge 1,\end{array} \right.
\]
which gives the cost of such a nucleus inside the domain.
We are instead interested in the case when $\Omega$ is a {\it generic corner}. Our main result reads as  follows:
\begin{theorem}
 If the domain $\Omega$ is a generic corner (see~\eqref{domain} for a precise definition), then the minimum of the energy scales like 
 \begin{equation*}
  \inf_{\chi \textrm{ satisfies } \eqref{chiBV},\eqref{volume}} E[\chi] \sim \left\{ \begin{array}{ll} V^{2/3} & \textrm{ if } V \le 1 \\ V^{7/9} & \textrm{ if } V \ge 1.\end{array} \right.
 \end{equation*}
\end{theorem}
This results show that it is indeed energetically favored to start nucleation at corners. We should however warn the reader that we do not say anything about the nucleation on faces or edges of the specimen. We believe that in these cases, the scaling should be $V^{9/11}$ as for the nucleation inside but since they do not seem accessible using the presently-available techniques
\footnote{The proof of the lower bound for the whole space in~\cite{KKO} is based on the Fourier transform, which for obvious reasons can not be used here.}, we do not investigate this question here. 
We would also like to point out that our result does not hold for the corner $(\R^3)^+$ since in this case there is no habit plane which is cutting the corner. By adapting the construction of \cite{KKO}, we see that in this case, nucleation both at the corner and at the edge can be achieved with an energy of order at most $\min(V^{2/3},V^{8/10})$. Since this case is very degenerate and since again, a proof of the corresponding lower bound would probably require completely new ideas, we do not push further this question here. \\

The main difference between nucleation at a corner and nucleation in the bulk is that in the latter, the notion of self-accommodation plays a central role.
 In fact, in order to create a nucleus in the bulk, one needs to be able to construct a three-dimensional macroscopically stress-free configuration.
 As shown in \cite{KKO,KO}, this is possible only if the three variants of the martensite are present (in similar volume fractions). Indeed, if only two variants of martensite are used, then the energy must be larger --  of order $V^{10/12}$) \cite{KO}.
 Moreover, to avoid large energetic costs the martensitic inclusion prefers to have boundary parallel to one of the habit planes. This is not completely possible 
for the inclusion in the whole space, and so one constructs a lens-shaped inclusion \cite{KKO} to (at least approximately) have the boundary parallel to 
the habit plane.   The energy scaling law $V^{9/11}$ is partly due to the {\it macroscopic bending} of the lens. 
On the contrary, for a generic corner it is possible to construct a macroscopically stress-free configuration with a boundary between martensitic inclusion and austenite exactly parallel to the habit plane. 
The optimal scaling of the energy $V^{7/9}$ is then obtained by the classical branching pattern using only two variants of martensite (see \cite{ KM:94,CapOt1,CKO}).   

Our proof of the lower bound strongly relies on a result by Capella and Otto \cite{CapOt2} which gives a local lower bound on the energy of austenite-martensite interfaces. More precisely, we use the fact that if in a ball the volume fraction of austenite is bounded away from $0$ and $1$ (i.e. both austenite and martensites are present), then in a larger ball the energy cannot be small. 

Let us emphasis that we are dealing only with the cubic-to-tetragonal phase transformation where the possible microstructures are very rigid (see \cite{CapOt1,CapOt2} and references therein). Much less is known in other situations such as cubic-to-orthorhombic or cubic-to-monoclinic phase transformations (see \cite{Angk1,Angk2} for some available results).
\\

\textbf{Notation} In the paper we will use the following notation. The symbols $\sim$, $\ges$, $\les$ indicate estimates that hold up to a constant, which could depend on the domain $\Omega$, but not on $\chi$ or $u$. For instance, $f\les g$ denotes the existence of a constant $C=C(\Omega)>0$ such that $f\le Cg$. For $u,v$ in $\R^3$, the tensor product $u\otimes v$ is the $3\times 3$
 matrix defined componentwise by $(u\otimes v)_{ij}:=u_iv_j$. The symmetrized tensor product is defined by $u\odot v:= \frac{1}{2}( u\otimes v+ v\otimes u)$. 
Finally, for a matrix $A$, we denote $e(A):=\frac{1}{2} (A^t+A)$ its symmetric part and  its norm by $|A|:=\sqrt{\sum_{ij} A_{ij}^2}$.  

\section{The model}\label{secmod}

We adopt the notation from \cite{KKO}. We assume that the atoms in the specimen are aligned with the coordinate axes, i.e.
\begin{equation*}
e^{(0)} := 0,\quad
 e^{(1)} := \begin{pmatrix} -2 & 0 & 0 \\ 0 & 1 & 0 \\ 0 & 0 & 1\end{pmatrix},\quad
 e^{(2)} := \begin{pmatrix}  1 & 0 & 0 \\ 0 &-2 & 0 \\ 0 & 0 & 1\end{pmatrix},\quad
 e^{(3)} := \begin{pmatrix}  1 & 0 & 0 \\ 0 & 1 & 0 \\ 0 & 0 &-2\end{pmatrix}.
\end{equation*}

Two strains $A,B$ will be called compatible if there exists a plane and a continuous function $u$ with $e(u)=A$ on one side of the plane and $e(u)=B$ on the other side. Two such strains are called twins, the corresponding plane is called the twinning plane and its normal is called the twin direction. It can be proven that two strains $A$ and $B$ are compatible if and only if
$A-B= a \odot b$ for some vectors $a, b$ in $\R^3$. In the case of cubic-to-tetragonal phase transformations it can be computed that the martensitic variants are pairwise compatible in the sense that for any permutation $(ijk)$, letting $\eps_{ijk}$ be the signature of the permutation, there holds
\[e^{(i)}-e^{(j)}=6 \eps_{ijk} (b_{ij} \odot b_{ji}),\]
where 
the six twinning planes have normals
\begin{equation*}
\begin{array}{lll}
 b_{12}:=\ds\frac{1}{\sqrt{2}}\begin{pmatrix}
                            1\\1\\0
                           \end{pmatrix}, &
b_{31}:=\ds\frac{1}{\sqrt{2}}\begin{pmatrix}
                            1\\0\\1
                           \end{pmatrix}, &
b_{23}:=\ds\frac{1}{\sqrt{2}}\begin{pmatrix}
                            0\\1\\1
                           \end{pmatrix}, \\[10pt]
b_{21}:=\ds\frac{1}{\sqrt{2}}\begin{pmatrix}
                            -1\\1\\0
                           \end{pmatrix}, &
b_{13}:=\ds\frac{1}{\sqrt{2}}\begin{pmatrix}
                            1\\0\\-1
                           \end{pmatrix}, &
b_{32}:=\ds\frac{1}{\sqrt{2}}\begin{pmatrix}
                            0\\-1\\1
                           \end{pmatrix}.
\end{array}
\end{equation*}
Moreover, though no single martensitic variant is compatible with the austenite, there is compatibility of austenite with some convex combinations of the martensite variants in the sense that for $i\neq j$,
\[\frac{1}{3} e^{(i)}+ \frac{2}{3} e^{(j)}= 2\eps_{ijk} ( b_{jk} \odot b_{kj}).\]
The planes which are connecting austenite to twinned martensite are called habit planes and their normals are called habit directions. Notice that for cubic-to-tetragonal phase transformation, habit planes and twinning planes coincide. We refer the interested reader to \cite{bhat} for more information on the subject. \\

To study the nucleation in the corner, we consider a domain $\Omega$ of the form
\begin{equation}\label{domain}
 \Omega := \left\{ \alpha a + \beta b + \gamma c : \alpha \ge 0,\beta\ \ge 0,\gamma \ge 0 \right\},
\end{equation}
where $a,b,c$ is some basis of $\R^3$ with $|a|=|b|=|c|=1$. Moreover, we assume that (at least) one of the habit planes can cut off a corner of our specimen, i.e. that at least one of the normals to the habit planes (we denote that normal by $n$) satisfies either
\begin{equation}\label{orientation}
 a \cdot n > 0, \quad b \cdot n > 0, \quad c \cdot n > 0, 
\end{equation}
or
\begin{equation*}
 a \cdot n < 0, \quad b \cdot n < 0, \quad c \cdot n < 0.
\end{equation*}
Without loss of generality we can assume that $n = b_{23}$ and that condition~\eqref{orientation} holds. Moreover, we can also assume that the vectors $a$, $b$, and $c$ are in positive order, i.e. that $(a \times b) \cdot c > 0$, $(b \times c) \cdot a > 0$, and $(c \times a) \cdot b > 0$. Finally, we define 
\[\mu := \min\left\{\frac{(a \times b)}{|a \times b|} \cdot c, \frac{(b \times c)}{|b \times c|} \cdot a, \frac{(c \times a)}{|c \times a|} \cdot b \right\} \in (0,1].\]


\section{Lower Bound}\label{seclow}

In this section we prove the lower bound, that is
\begin{equation}\label{lb}
 \inf_{\chi \textrm{ satisfies } \eqref{chiBV},\eqref{volume}} E[\chi] \gtrsim \left\{ \begin{matrix} V^{2/3} & \textrm{ if } V \le 1, \\ V^{7/9} & \textrm{ if } V \ge 1.\end{matrix} \right.
\end{equation}
If $V \le 1$, the lower bound immediately follows from the isoperimetric inequality. Hence we need to focus only on the case $V \ge 1$. In the proof we will use the following lemma, which can be obtained by simple scaling from \cite[Theorem 1 - part i)]{CapOt2}:
\begin{lemma}\label{lm1}
There exists a small but universal $\kappa \in (0,1)$ with the following property: for every $R \ge 1/\kappa$, every $u \in H^1(B_R(0),\R^3)$ and every $\chi$ with the property
 \begin{equation}\label{mixedball}
  \frac{1}{10} |B_{\kappa R}| \le \int_{B_{\kappa R}} \chi_1 + \chi_2 + \chi_3 \dx \le \frac{9}{10} |B_{\kappa R}|,
 \end{equation}
we have 
\begin{equation}\label{lbball}
 \int_{B_R} \left|e(u) - \sum_{i=1}^3 \chi_i e^{(i)}\right|^2 + \sum_{i=1}^3 |\nabla \chi_i| \dx \gtrsim R^{7/3}.
\end{equation}
\end{lemma}

It is easy to see that for $R < 1/\kappa$, \eqref{lbball} (with possibly different universal constant) follows from~\eqref{mixedball} by the isoperimetric inequality. 

Let us introduce the following notation for the portion of the space occupied by the  martensites
\begin{equation}
 M := \left\{ x \in \Omega : \chi_1(x) + \chi_2(x) + \chi_3(x) = 1 \right\},
\end{equation}
so that $|M|=V$. Let us now consider $\chi$ which satisfies~\eqref{chiBV},~\eqref{volume}. To prove~\eqref{lb} it is enough to construct a (at most countable) set of balls $B(x_i,r_i)$ with the following properties:
\begin{enumerate}
 \item the balls $B(x_i,r_i)$ are disjoint and are subset of $\Omega$; 
 \item $\left| M \setminus \bigcup_{i} B(x_i,15\mu^{-1}r_i) \right| = 0$;
 \item 
  $\frac{1}{10} |B(x_i,\kappa r_i)| \le \int_{B(x_i,\kappa r_i)} \chi_1 + \chi_2 + \chi_3 \dx \le \frac{9}{10} |B(x_i,\kappa r_i)|.$
\end{enumerate}
Indeed, let us assume that we have such a covering. Then by Lemma~\ref{lm1} we have for every $i$
\begin{equation*}
 \int_{B(x_i,r_i)} \left|e(u) - \sum_{i=1}^3 \chi_i e^{(i)}\right|^2 + \sum_{i=1}^3 |\nabla \chi_i| \dx \gtrsim r_i^{7/3}.
\end{equation*}
Since $B(x_i,r_i)$ are disjoint and $B(x_i,r_i) \subset \Omega$, summing the above relation in $i$ implies
\begin{equation*}
 \int_\Omega \left|e(u) - \sum_{i=1}^3 \chi_i e^{(i)}\right|^2 + \sum_{i=1}^3 |\nabla \chi_i| \dx \gtrsim \sum_{i} r_i^{7/3}.
\end{equation*}
Finally, because of the second property, we see that 
\begin{equation*}
 \sum_i r_i^3 \gtrsim \sum_i |B(x_i,15\mu^{-1}r_i)| \ge \left| \bigcup_i B(x_i,15\mu^{-1}r_i)\right| \ge |M| = V,
\end{equation*}
and by Jensen's inequality applied to the concave function $t^{7/9}$, we see that
\begin{equation*}
 \int_\Omega \left|e(u) - \sum_{i=1}^3 \chi_i e^{(i)}\right|^2 + \sum_{i=1}^3 |\nabla \chi_i| \dx \gtrsim \sum_{i} (r_i^3)^{7/9} \gtrsim \left( \sum_{i} r_i^3 \right)^{7/9} \gtrsim V^{7/9}.
\end{equation*}
Since this holds for any $u$ and $\chi$ which satisfies~\eqref{chiBV} and~\eqref{volume}, \eqref{lb} immediately follows.

It remains to show the existence of the balls $B(x_i,r_i)$. Let $\bar x$ be any point of density of $M$.\footnote{That is a Lebesgue point of $\chi_M$. Notice that by definition, $\bar x$ has to belong to the interior of $\Omega$.} Then there exists a radius $r_0 > 0$ such that $B(\bar x,r_0) \subset \Omega$ and $|B(\bar x,\kappa r_0) \cap M| \ge \frac{1}{10} |B(\bar x,\kappa r_0)|$. Consider the function $f$ defined by
\begin{equation*}
f(t) := \frac{|B(\bar x + t(a+b+c),\kappa(r_0 + \mu t)) \cap M|}{|B(\bar x + t(a+b+c),\kappa(r_0 + \mu t))|},
\end{equation*}
which measures the volume fraction of balls obtained from $B(\bar x,r_0)$ by dilation in such a way that all these balls belong to $\Omega$. 

We claim that for some $t_{\bar x} \ge 0$ we have $\frac{1}{10} \le f(t_{\bar x}) \le \frac{9}{10}$. Indeed, we know that $f(0) \ge \frac{1}{10}$ and that $f$ is a continuous function on $[0,\infty)$. Since the volume of $M$ is finite, we have $\lim_{t \to \infty} f(t) = 0$, and the claim follows. In this way, to each $\bar x$ we assign a ball $B(\bar x + t_{\bar x}(a+b+c),r_0 + \mu t_{\bar x})$. Before proceeding we point out that $B(\bar x + t_{\bar x}(a+b+c),r_0 + \mu t_{\bar x}) \subset \Omega$. 

Now let $\mathcal{B}$ denotes set of the constructed balls enlarged by a factor $3\mu^{-1}$,  i.e. we define $\mathcal{B} := \left\{ B(\bar x + t_{\bar x}(a+b+c),3\mu^{-1}(r_0 + \mu t_{\bar x})) : \bar x \textrm{ point of density of } M\right\}$. We note that $\bar x \in B(\bar x+t_{\bar x}(a+b+c), 3\mu^{-1} (r_0 + \mu t_{\bar x}))$, in particular $M$ is covered (up to a negligible set) by the balls in $\mathcal{B}$. 

Now we apply Vitali's covering lemma to $\mathcal{B}$. We observe that for every $\bar x$ and corresponding $t_{\bar x}$ we have $|M| \ge \left|B(\bar x + t_{\bar x}(a+b+c),\kappa(r_0 + \mu t_{\bar x})) \cap M\right| \ge \frac{1}{10} 
|B(\bar x + t_{\bar x}(a+b+c),\kappa(r_0 + \mu t_{\bar x}))|$, and so $r_0 + \mu t_{\bar x} \lesssim |M|^{1/3}$. Therefore the supremum of radii of balls in $\mathcal{B}$ is bounded, and so by Vitali's covering lemma there exists a set of balls $B_i = B(x_i,r_i) = B(x_i + t_{x_i}(a+b+c),r_0 + \mu t_{x_i}) \subset \Omega$ such that $B(x_i,3\mu^{-1}r_i) \in \mathcal{B}$, $B(x_i,15 \mu^{-1}r_i)$ cover $M$ (up to a negligible set), and $B(x_i,3 \mu^{-1} r_i)$ (hence also $B(x_i,r_i)$) are disjoint. This completes the proof of the lower bound. 

\section{Upper bound}\label{secup}
In this section we prove the upper bound, i.e. 
\begin{equation}\label{ub}
 \inf_{\chi \textrm{ satisfies } \eqref{chiBV},\eqref{volume}} E[\chi] \les \left\{ \begin{matrix} V^{2/3} & \textrm{ if } V \le 1, \\ V^{7/9} & \textrm{ if } V \ge 1.\end{matrix} \right.
\end{equation}
If $V\le1$, it is enough considering $\chi_1=\chi_{B_R(x)}$ with $|B_R(x)|=V$ and $B_R(x)\subset \Omega$, $\chi_2=\chi_3=0$ and $u=0$ to get 
\[E(\chi)\le CV^{2/3}+ \int_{B_R(x)} |e^{(1)}|^2\les V^{2/3}+V\les V^{2/3}.\]
For $V\ge 1$, we use a classical branching construction which dates back to Kohn and M\"uller \cite{KM:94} (see also \cite{CapOt1,CKO}). Our construction will be essentially two-dimensional and will use branching close to the martensite-austenite interface. We will use a simplified version of the construction used in~\cite{KKO} and will follow their notations. \\

Let $n:=b_{32}$ be the normal to the habit plane cutting off the corner and let us denote 
\begin{equation*}
 b_3:=\frac{b_{21}\times n}{|b_{21}\times n|}, \qquad b_2:=\frac{n\times b_3}{|n\times b_3|}, \qquad b_1:=\frac{b_3\times b_{21}}{|b_3\times b_{21}|},
\end{equation*}
 which is a basis of $\R^3$ and let $y_i:= x\cdot b_i$ be the associated coordinates. For $R>1$, let $C_R:=\{-R\le y_1\le R, 0\le y_2\le R, 0\le y_3\le R\}$.  
We are going to construct a pair $(u,\chi)$ such that $\chi_1 +\chi_2+\chi_3=1$ in $C_R^-:=C_R\cap\{y_1\le 0\}$, $\chi_1+\chi_2+\chi_3=0$ in $C_R^+:=C_R\cap\{y_1> 0\}$, $u=0$ in $C_R^+$ and
\begin{equation}\label{estimsupfonda}
 \int_{C_R}  \left( |\nabla \chi_1| + |\nabla \chi_2| + |\nabla \chi_3| \right) \dx+ \int_{C_R} \left| e(u) - \sum_{i=1}^3 \chi_i e^{(i)} \right|^2 \dx \les R^{7/3}.
\end{equation}
 Let us first see how this would give the proof of \eqref{ub}. By the non-degeneracy condition \eqref{orientation}, we see that for every $V>1$ we can find $R\sim V^{1/3}$ and a translation $z_R\in \R^3$ such that $|\Omega \cap (z_R+C_R^-)|=V$ and $\{ x \in \Omega : x \cdot n < 0 \} \subset z_R + C_R^-$. Then by extending the previously constructed pair $(u,\chi)$ to the whole $\Omega$ by zero we find that 
\[E(\chi)\le  \int_{C_R}  \left( |\nabla \chi_1| + |\nabla \chi_2| + |\nabla \chi_3| \right) \dx+ \int_{C_R} \left| e(u) - \sum_{i=1}^3 \chi_i e^{(i)} \right|^2 \dx \les R^{7/3}\les V^{7/9}\]
and the upper bound is proven. Let us now turn to the construction of the pair $(u,\chi)$ satisfying \eqref{estimsupfonda}. In our construction we are going to use only variants $1$ and $2$ of the martensite, i.e. $\chi_3\equiv 0$. 
In the martensitic phase $C_R^-$, we will have fine-scale oscillations of the martensite variants in the direction $b_2$  and branching in the direction $b_1$. The whole construction will be invariant in the direction $b_3$. 
As in \cite{KKO} first we need to choose the gradients that will be involved and which allow for twinning between martensites and compatibility with the austenite. For this we let
\[D^{(1)}:=\begin{pmatrix}
              -2 & 2 & 0\\
	    -2& 1 & 1\\
	      0&-1 & 1
             \end{pmatrix}
\quad D^{(2)}:=\begin{pmatrix}
              1 & -1 & 0\\
	    1& -2 & 1\\
	      0&-1 & 1
             \end{pmatrix}
  \quad D^{M}:=\begin{pmatrix}
              0 & 0 & 0\\
	    0& -1 & 1\\
	      0&-1 & 1
             \end{pmatrix},
\]
so that $e(D^{(1)})=e^{(1)}$, $e(D^{(2)})=e^{(2)}$, $D^{(1)}-D^{(2)}= 6(b_{12}\otimes b_{21})$, and 
\[
 D^M=\frac{1}{3} D^{(1)}+\frac{2}{3} D^{(2)}=2 (b_{23}\otimes n).
\]
In $C_R^-$ we are going to let $u=u^m+D^M$, where $u^m$ will be the microscopic displacement accounting for the twinning of the martensites and for the branching process, and where $D^M$ is the macroscopic displacement, which ensures compatibility with the austenite phase. 

Since the construction consists of self-similar cells, let us now define the basic cell $Z$. For $1< h\sim w^{3/2}$ (which would then imply $w\les h$), let
\[Z:=\{ 0\le y_1\le h, 0\le y_2\le w, 0\le y_3 \le w\}.\]
We let $\chi_1:=1$ on the sets
\[\left\{ \left|\frac{y_2}{w}-\frac{1}{6}\right|\le \frac{y_1}{18 h}\right\}\cap Z, \quad  \left\{ \left|\frac{y_2}{w}-\frac{5}{6}\right|\le \frac{y_1}{18 h}\right\}\cap Z, \quad 
 \left\{ \left|\frac{y_2}{w}-\frac{1}{2}\right|\le \frac{1}{6}-\frac{y_1}{9 h}\right\}\cap Z, 
\]
and $\chi_1 :=0$ on the rest of $Z$. We then let $\chi_2=1-\chi_1$. Notice that on each slice $\{y_1=c\}$, there holds
\[ \int_{\{y_1=c\} \cap Z} \chi_1 =\frac{w^2}{3} \qquad \textrm{and} \qquad \int_{\{y_1=c\} \cap Z} \chi_2 =\frac{2w^2}{3}. \]
We define the microscopic displacement $u^m$ by imposing that 
\begin{equation}\label{umvanishes}
 u^m = 0 \ \textrm{ on } \ \{y_2 \in \{0,w\} \textrm{ or } y_3 \in\{0,w\}\} \cap Z
\end{equation}
and so that the derivatives of $u^m$ in the $b_2$ and $b_3$ directions are given by
\begin{align*}
 \partial_{b_2} u^m&:= [(D^{(1)}- D^M)\chi_1 +(D^{(2)}-D^M) \chi_2] b_2,\\
\partial_{b_3} u^m&:= [(D^{(1)}- D^M)\chi_1 +(D^{(2)}-D^M) \chi_2] b_3=0.
\end{align*}
This together with~\eqref{umvanishes} implies that $\partial_{b_1} u^m$ is constant on each connected component of the support of $\chi_1$ and $\chi_2$, and has a jump of order $\frac{w}{h}$ at the interfaces. Let us now estimate the energy of such a configuration. Since $(D^{(1)}- D^M)b_1=(D^{(2)}- D^M)b_1=0$, the definition of $\partial_{b_2} u $ and 
$\partial_{b_3} u$ imply
\[\int_{Z} |e(u^m)-\sum_{i=1}^2 \chi_i(e^{(i)}-e(D^M)|^2 \le \int_{Z} |\nabla u^m-\sum_{i=1}^2 \chi_i(D^{(i)}-D^M|^2=\int_Z |\partial_{b_1} u|^2 \les \frac{w^4}{h}.\] 
Since $w\les h$, the interfacial energy can be estimated by $w h$ and thus
\begin{equation}\label{estimenercell}
\int_Z \sum_{i=1}^2 |\nabla \chi_i| +\int_{Z} |e(u^m)-\sum_{i=1}^2 \chi_i(e^{(i)}-e(D^M)|^2 \les h w+\frac{w^4}{h} \sim w^{5/2}.
\end{equation}

We now decompose $C_R^-$ into cells. Let us first denote by $C_R^{bl}:=C_R^-\cap\{-1\le y_1\le 0\}$ the boundary layer of thickness $1$ and by $C_R^{int}:= C_R^-\cap\{-R\le y_1\le -1\}$ the interior domain. We then decompose the square $\{ 0\le y_2\le R, 0\le y_3\le R\}$ into cubes $Q_k$, $k=1,.., K$ of sidelength $w\sim R^{2/3}$ (so that $K\sim R^{2/3}$) and consider the corresponding cylinder
\[\Sigma_k:=\left\{ q+\alpha b_1 : q\in Q_k, -R\le \alpha \le -1\right\}\subset C_R^-. \]

Next, we decompose any such cylinder in refining cells which are going to be rescaled versions of $Z$. The first cell is of width $w$ and height $h_0:=C_1 w^{3/2}$ (with $C_1$ to be fixed later) and the i-th generation of cells is defined by
\[w_i:=\frac{w_{i-1}}{3} \qquad \textrm{ and } \qquad h_i:=C_1 w_i^{3/2}\]  
so that above each cell of the generation $i-1$ there are nine cells of the generation $i$. The construction stops after $M$ iterations when 
\begin{equation}\label{stopcrit}
 h_M\le w_M.
\end{equation}
The constant $C_1$ is finally chosen so that $\sum_{i=0}^{M} h_i =R-1$. Notice that since $h_i$ is geometric and since $w\sim R^{2/3}$, we have $C_1\sim 1$. The functions $\chi_1$, $\chi_2$ and $u^m$ are then defined on the constructed cells by rescaling of those defined in $Z$. 
In the boundary layer, we set $\chi_1=1$ and extend continuously $u^m$ so that $u=0$ on $\{y_1=0\}$ and
\[\|\nabla u^m\|_{L^{\infty}(C_R^{bl})}\les \|u^m\|_{L^{\infty}( \partial C_R^{bl})}+ \|\nabla u^m\|_{ L^{\infty}( \partial C_R^{bl})}. \]
We can now compute the energy of such a configuration. The energy can be split into two parts, one coming from the contribution in $C_R^{int}$ and the other coming from the contribution in $C_R^{bl}$. Let us start by estimating the energy coming from $C_R^{int}$. For this, consider first one cylinder $\Sigma_k$. Recalling \eqref{estimenercell} and the definitions of $w_i$ and $h_i$, we find
\[\int_{\Sigma_k} \sum_{i=1}^2 |\nabla \chi_i| +\int_{\Sigma_k} |e(u^m)-\sum_{i=1}^2 \chi_i(e^{(i)}-e(D^M))|^2 \les w^{5/2}\sum_{i=0}^{+\infty} 3^{2i} \left(\frac{1}{3}\right)^{5i/2}\les R^{5/3}\]
which, summing for $k=1,..,K$ (and recalling that $K\sim R^{2/3}$) gives
\[\int_{C_R^{int}} \sum_{i=1}^2 |\nabla \chi_i| +\int_{C_R^{int}} |e(u^m)-\sum_{i=1}^2 \chi_i(e^{(i)}-e(D^M))|^2 \les R^{7/3}.\]

We finally estimate the contribution of the energy coming from the boundary layer. Notice that since its thickness is one, $|C_R^{bl}|\sim R^2$ and $\mathcal{H}^{2}(\partial C_R^{bl})\sim R^2$. Moreover, since in the last generation of cells $w_M\sim 1$, $\|u^m\|_{L^{\infty}( \partial C_R^{bl})}+ \|\nabla u^m\|_{ L^{\infty}( \partial C_R^{bl})}\les 1$ and therefore
\[\int_{C_R^{int}} \sum_{i=1}^2 |\nabla \chi_i| +\int_{C_R^{int}} |e(u^m)-\sum_{i=1}^2 \chi_i(e^{(i)}-e(D^M))|^2 \les R^2\les R^{7/3},\]
from which \eqref{estimsupfonda} follows.
\section*{Acknowledgment}
The authors warmly thank F.Otto for suggesting the problem and for very valuable discussions. They also thank H. Seiner for many comments on the experimental and physical background. M. Goldman was funded by the Von Humboldt foundation. 
\bibliography{corner}

\end{document}